\documentclass[11pt,twoside]{article}
\usepackage[plainpages=false]{hyperref}
\usepackage{amsfonts,latexsym,rawfonts,amsmath,amssymb,amsthm}
\usepackage{amsmath,amssymb,amsfonts,latexsym,lscape,rawfonts,authblk}
\numberwithin{equation}{section}

\newcommand{\pd}[2]{\frac {\partial #1}{\partial #2}}

\newcommand{\oo}{\omega}

\newcommand{\Na}{\nabla}

\newcommand{\ee}{\epsilon}

\newcommand{\beq}{\begin{equation}}
\newcommand{\eeq}{\end{equation}}
\newcommand{\beqs}{\begin{eqnarray*}}
\newcommand{\eeqs}{\end{eqnarray*}}
\newcommand{\beqn}{\begin{eqnarray}}
\newcommand{\eeqn}{\end{eqnarray}}
\newcommand{\beqa}{\begin{array}}
\newcommand{\eeqa}{\end{array}}

\def\lra{\longrightarrow}

\def\bc{\begin{center}}
\def\ec{\end{center}}

\def\i{\sqrt{-1}}

\def\cP{{\cal P}}
\def\cW{{\cal W}}

\def\cP{{\mathcal P}}

\def\cW{{\mathcal W}}

\def\RR{{\mathbb R}}

\def\CC{{\mathbb C}}

\def\begeq{\begin{equation}}
\def\endeq{\end{equation}}
\def\and{\quad{\rm and}\quad}

\let\lra=\longrightarrow

\def\mapright\#1{\,\smash{\mathop{\lra}\limits^{\#1}}\,}

\def\an{\;\;\;{\rm and}\;\;\;}

\def\ri{\rightarrow}

\def\un{\underline}

\def\pbp{\sqrt{-1}\partial\bar\partial}

\newtheorem{prop}{Proposition}[section]
\newtheorem{theo}[prop]{Theorem}
\newtheorem{lem}[prop]{Lemma}

\title{{\bf\Large{On the lower bound of the $K$-energy and $F$-functional\footnote{
After we posted our first version of this paper on the arXiv, we got
some feedbacks that our result can be proved  by the continuity
method. However, the idea of the proof, which comes from our joint
paper \cite{[CLW]}, is still interesting, and may have some other
applications.}}}}
\author{Haozhao Li }
\date{}

\begin{document}
\bibliographystyle{plain}
\maketitle \footnotetext{{Math Subject Classifications:}
           32Q20, 53C44.}

\begin{abstract}Using Perelman's results on K\"ahler-Ricci flow, we  prove that the $K$-energy is bounded from below
if and only if the $F$-functional is bounded from below in the
canonical K\"ahler class.
\end{abstract}

\section{Introduction}
One of the central problems in K\"ahler geometry is to study the
existence of K\"ahler-Einstein metrics, which is closely related to
the behavior of several energy functionals. During the last few
decades, these energy functionals  have been intensely studied and
there are many interesting results. The $K$-energy, which was
introduced by Mabuchi in \cite{[Ma]}, plays an important role in
K\"ahler geometry.

Let $(M, \oo)$ be an $n$-dimensional compact K\"ahler manifold with
$c_1(M)>0$. We define the space of K\"ahler potentials by
$$\cP(M, \oo)=\{\varphi\in C^{\infty}(M, \RR)\;|\; \oo+\pbp\varphi>0\},$$
where $\oo\in 2\pi c_1(M)$. For any $\varphi\in \cP(M, \oo)$, we
define the $K$-energy by \beq \nu_{\oo}(\varphi)=-\frac 1V
\int_0^1\;\int_M\; \pd {\varphi_t}{t}(R_{\varphi_t}-\un
R)\oo_{\varphi_t}^n\wedge dt\label{eq:K}\eeq where $\varphi_t(t\in
[0, 1])$ is a path in $\cP(M, \oo)$ with $\varphi_0=0$ and
$\varphi_1=\varphi$,  $\un R$ is the average of scalar curvature,
and $V=[\oo]^n$ is the volume. Bando-Mabuchi \cite{[BaMa]} showed
that if $M$ admits a K\"ahler-Einstein metric, then the $K$-energy
is bouned from below. Later, Tian \cite{[Tian97]}\cite{[Tianbook]}
proved that the existence of K\"ahler-Einstein metrics is equivalent
to the properness of the $K$-energy in the canonical K\"ahler class.
In fact, Tian proved that the existence of K\"ahler-Einstein metrics
is equivalent to the properness of the $F$-functional, which was
introduced by Ding-Tian \cite{[DingTian]} as follows \beq
F_{\oo}(\varphi)=\frac 1V\sum_{i=0}^{n-1}\; \frac {i+1}{n+1}\int_M\;
\i\partial \varphi\wedge \bar\partial \varphi\wedge \oo^i\wedge
\oo_{\varphi}^{n-1-i}-\frac 1V\int_M\; \varphi \oo^n-\log \Big(\frac
1V\int_M\; e^{h_{\oo}-\varphi}\oo^n\Big).\label{eq:F}\eeq

To prove the convergence of K\"ahler-Ricci flow, Chen-Tian
\cite{[ChenTian]}\cite{[chen-tian1]} introduced a series of energy
functionals $E_k(k=0, 1, \cdots, n)$ defined by \beqs E_{k,
\oo}(\varphi)&=&\frac 1V\int_M\; \Big(\log \frac
{\oo^n_{\varphi}}{\oo^n}-h_{\oo}\Big)\Big( \sum_{i=0}^k\;
Ric_{\varphi}^i\wedge \oo^{k-i}\Big)\wedge \oo_{\varphi}^{n-k}+\frac
1V\int_M\; h_{\oo} \Big(\sum_{i=0}^k\; Ric_{\oo}^i\wedge
\oo^{k-i}\Big)\wedge \oo^{n-k}\\&&+\frac {n-k}{V}\int_0^1\;\int_M\;
\pd {\varphi_t}{t}(\oo_{\varphi_t}^{k+1}-\oo^{k+1})\wedge
\oo_{\varphi_t}^{n-k-1}\wedge dt,\eeqs where  $h_{\oo}$ is the Ricci
potential defined by \beq Ric(\oo)-\oo=\pbp h_{\oo}, \an
\int_M\;(e^{h_{\oo}}-1)\oo^n=0,\label{eq0.1}\eeq and $\varphi_t(t\in
[0, 1])$ is a path from $0$ to $\varphi$ in $\cP(M, \oo)$. The first
energy $E_0$ of these series is exactly the $K$-energy, and the
second $E_1$ is the Liouville energy on Riemann surfaces.

There are many relations between these energy functionals. Pali
\cite{[Pali]} prove that $E_1$ is bounded from below if the
$K$-energy is bounded from below. Recently, Chen-Li-Wang
\cite{[CLW]} proved the converse is also true. There are also some
results on the lower bound of $E_k.$ Following a question proposed
by X. X. Chen \cite{[Chen1]}, Song-Weinkove   \cite{[SoWe]} showed
that the existence of K\"ahler-Einstein metrics is equivalent to the
properness of $E_1$ in the canonical class, and they also showed
that $E_k$ are bounded from below under some additional curvature
conditions. Recently, following suggestion of X. X. Chen, the author
\cite{[Li]}  found new relations between all these functionals and
generalized Pali-Song-Weinkove's results.

In summary, the relations between the existence of K\"ahler-Einstein
metrics and these energy functionals can be roughly written as
follows:  $M$ admits K\"ahler-Einstein metrics $\Longleftrightarrow$
the $F$-functional is proper $\Longleftrightarrow$ the $K$-energy is
proper $\Longleftrightarrow$ $E_1$ is proper. A natural question is
what will happen if these energy functionals are just bounded from
below instead of proper. In this paper, we prove

\begin{theo}\label{main}The $K$-energy is
bounded from below if and only if $F$ is bounded from below on
$\cP(M, \oo).$ Moreover, we have
$$\inf_{\oo'\in [\oo]}F_{\oo}(\oo')=\inf_{\oo'\in [\oo]}\nu_{\oo}(\oo')-\frac 1V\int_M\; h_{\oo}\oo^n,$$
where $h_{\oo}$ is the Ricci potential with respect to the metric
$\oo.$
\end{theo}
Combining this with the results in \cite{[CLW]}, we actually prove
that $F$ is bounded from below $\Longleftrightarrow$ the $K$-energy
is bounded from below $\Longleftrightarrow$ $E_1$ is bounded from
below. We expect that the lower boundedness of all energy
functionals $E_k$ is equivalent, and perhaps the lower boundedness
implies the existence of singular K\"ahler-Einstein metrics and
certain stabilities. \\

The idea of the proof of Theorem \ref{main} is essentially due to
our joint paper \cite{[CLW]}. The key point is to  estimate the
difference of $F$ and $\nu_{\oo}$ along the K\"ahler-Ricci flow, and
we show that the difference of these two functionals at infinity is
a uniform constant independent of the initial metric of the flow.
However, the proof needs Perelman's deep estimates on the
K\"ahler-Ricci flow, while in \cite{[CLW]} the equivalence of the
$K$-energy and $E_1$ doesn't. This is because we can compare the
derivatives of these energy functionals along the K\"ahler-Ricci
flow in \cite{[CLW]}, but we don't have similar estimates in this
paper. The readers are referred to \cite{[CLW]} for details. We
expect that this flow method can be used to
prove the equivalence of all $E_k$ functionals in the future.\\

\noindent {\bf Acknowledgements}:  I would like to thank my advisors
Professor Xiuxiong Chen and Weiyue Ding for  their help and
encouragement over the past few years. I would also like to thank
Professor Gang Tian and Xiaohua Zhu for their help and some
enlightening discussions.

\section{K\"ahler-Ricci flow and the $K$-energy}
Let $(M, \oo)$ be an $n$-dimensional compact K\"ahler manifold with
$\oo\in 2\pi c_1(M)>0$. The  K\"ahler-Ricci flow with the initial
metric $\oo_0=\oo+\pbp \varphi_0$ is of the form
\begin{equation}
{{\partial \,\oo_{\varphi}} \over {\partial t }} = \oo_{\varphi} -
Ric_{\varphi}, \qquad \varphi(0)=\varphi_0.
\label{eq:kahlerricciflow}
\end{equation}
It follows that on the level of K\"ahler potentials, the
K\"ahler-Ricci flow becomes
\begin{equation}
{{\partial \varphi} \over {\partial t }} =  \log
{{\omega^n_{\varphi}} \over {\omega}^n } + \varphi - h_{\omega} ,
\label{eq:flowpotential}
\end{equation}
where $h_{\omega}$ is defined by (\ref{eq0.1}). Notice that for any
solution $\varphi(t)$ of (\ref{eq:flowpotential}), the function
$\tilde\varphi(t)=\varphi(t)+Ce^t$ is also a solution for any
constant $C$. Since
$$\pd {\tilde \varphi}{t}(0)=\pd {\varphi}{t}(0)+C,$$ we have
$$\frac 1V\int_M\; \pd {\tilde \varphi}{t}\oo_{\tilde
\varphi}^n\Big|_{t=0}=\frac 1V\int_M\; \pd
{\varphi}{t}\oo_{\varphi}^n\Big|_{t=0}+C.$$ Thus we can normalize
the solution $\varphi(t)$ such that the average of $\pd
{\varphi}t(0)$ is any
given constant.\\

Next we recall some basic facts on energy functionals.  The
$K$-energy, which is defined by (\ref{eq:K}), can be explicitly
expressed as (cf. \cite{[Chen0]}\cite{[Tianbook]}) \beq
\nu_{\oo}(\varphi)=\frac
1V\int_M\;\log\frac{\oo_{\varphi}^n}{\oo^n}\oo_{\varphi}^n+\frac
1V\int_M\; h_\oo(\oo^n-\oo_{\varphi}^n)-\frac
1V\sum_{i=0}^{n-1}\frac {n-i}{n+1}\int_M\; \i\partial \varphi\wedge
\bar\partial \varphi\wedge \oo^i\wedge
\oo_{\varphi}^{n-1-i}.\label{eq:eK}\eeq By direct calculation, the
$K$-energy is decreasing along the K\"ahler-Ricci flow. In fact, for
the solution $\varphi(t)$ of (\ref{eq:flowpotential}) we
have\footnote{Throughout this paper, the expressions such as $|\Na
f|$ and  $\Delta f$ are with respect to the metric
$\oo_{\varphi(t)}$. } \beq \frac {d}{dt}\nu_{\oo}(\varphi(t))=-\frac
1V\int_M\; \Big|\Na \pd {\varphi}t\Big|^2\oo_{\varphi}^n\leq 0.
\label{eq:dK}\eeq

 The following lemma tells us that if the $K$-energy is bounded from
below, we can normalize the solution such that the average of $\pd
{\varphi}{t}$ can be controlled. Since the normalization is crucial
in section 3, we include a proof here.

\begin{lem}\label{lem2.6}(cf. \cite{[ChenTian]})Suppose that the $K$-energy is bounded from below along the K\"ahler-Ricci flow. Then we can normalize the solution $\varphi(t)$ so that
$$c(0)=\frac 1V\int_0^{\infty}\;e^{-t}\int_M\;\Big|\Na \pd {\varphi}t\Big|^2\oo_{\varphi}^n\wedge dt<\infty, $$
where $c(t)=\frac 1V\int_M\; \pd {\varphi}t\oo_{\varphi}^n$. Then
for all time $t>0$, we have
$$c(t)>0,\quad \int_0^{\infty}\;c(t)dt<\nu_{\oo}(0)-\nu_{\oo}(\infty),$$
where $\nu_{\oo}(\infty)=\lim_{t\ri \infty}\nu_{\oo}(t)$.
\end{lem}
\begin{proof} A simple calculation yields
$$c'(t)=c(t)-\frac 1V\int_M\;|\Na \dot\varphi|^2\oo_{\varphi}^n.$$
Define
$$\ee(t)=\frac 1V\int_M\;|\Na \dot\varphi|^2\oo_{\varphi}^n.$$
Since the $K$ energy has a lower bound along the K\"ahler Ricci
flow, we have
$$\int_0^{\infty}\;\ee(t)dt=\frac 1V\int_0^{\infty}\int_M\;|\Na \dot\varphi|^2\oo_{\varphi}^n\wedge
dt= \nu_{\oo}(0)-\nu_{\oo}(\infty).$$ Now we normalize our initial
value of $c(t)$ as
\beqs c(0)&=&\int_0^{\infty}\;\ee(t)e^{-t}dt\\
&=&\frac 1V\int_0^{\infty}\;e^{-t}\int_M\;|\Na \dot\varphi|^2\oo_{\varphi}^n\wedge dt\\
&\leq &\frac 1V\int_0^{\infty}\int_M\;|\Na \dot\varphi|^2\oo_{\varphi}^n\wedge dt\\
&= &\nu_{\oo}(0)-\nu_{\oo}(\infty). \eeqs  From the equation for
$c(t)$, we have
$$(e^{-t}c(t))'=-\ee(t)e^{-t}.$$
Thus, we have \beqs0<c(t)=\int^{\infty}_t
\;\ee(\tau)e^{-(\tau-t)}d\tau \leq \nu_{\oo}(0)-\nu_{\oo}(\infty)
\eeqs and
$$\lim_{t\ri \infty}c(t)=\lim_{t\ri \infty}\int^{\infty}_t \;\ee(\tau)e^{-(\tau-t)}d\tau=0.$$
Since the $K$ energy is bounded from below, we have
$$\int_0^{\infty}\;c(t)dt=\frac 1V\int_0^{\infty}\int_M\;|\Na \dot\varphi|^2\oo_{\varphi}^n\wedge dt-c(0)\leq
\nu_{\oo}(0)-\nu_{\oo}(\infty).$$
\end{proof}

Now we recall the following result, which was proved by Perelman
using the $\cW$-functional and the gradient estimates for $\pd
{\varphi}t.$
\begin{lem}\label{Pere}(cf. \cite{[Natasa]}\cite{[NaTi]})For the solution $\varphi(t)$ of (\ref{eq:flowpotential}),
we choose $a_t$ by the condition $h_t=-\pd {\varphi}t+a_t$ such that
\beq\int_M\; e^{h_t} \oo_{\varphi}^n=V.\label{eq:P1}\eeq Then there
is a uniform constant $A$ independent of $t$ such that \beq
|h_t|\leq A,\;\;\;|\Na h_t|^2(t)\leq A,\an |\Delta h_t|\leq
A.\label{eq:P2}\eeq
\end{lem}

Finally, we state the following Poincar\'e inequality, which is
well-known in literature (cf. \cite{[Futaki]},\cite{[TianZhu]}).
\begin{lem}\label{TZ} For any K\"ahler metric $\oo_g$ and any
function $\phi\in C^{\infty}(M, \CC)$, we have
$$\int_M\; |\Na \phi|^2e^h\oo^n_g\geq \int_M\; |\phi-\un\phi|^2e^h \oo^n_g,$$
where $h$ is the Ricci potential function with respect to $\oo_g$
and
$$\un \phi=\frac 1V\int_M\; \phi e^h\oo_g^n.$$
\end{lem}
\vskip 1cm

\section{Proof of Theorem \ref{main}}
In this section, we prove the main theorem. First, by the expression
(\ref{eq:eK}) and (\ref{eq:F}), we can show the following lemma,
which directly implies the $K$-energy is bounded from below if $F$
is bounded from below.

\begin{lem}\label{lem3.1}\beq \nu_{\oo}(\varphi)\geq F_{\oo}(\varphi)+\frac
1V\int_M\; h_{\oo}\oo^n. \label{eq:3.1}\eeq
\end{lem}
\begin{proof}By the expression (\ref{eq:eK}),  the $K$-energy can be written as
 \beq\nu_{\oo}(\varphi)=\frac 1V\int_M\; u\oo_{\varphi}^n-\frac
1V\int_M\;\varphi\oo_{\varphi}^n+\frac 1V\int_M\;h_{\oo}\oo^n-\frac
1V\sum_{i=0}^{n-1}\frac {n-i}{n+1}\int_M\; \i\partial \varphi\wedge
\bar\partial \varphi\wedge \oo^i\wedge \oo_{\varphi}^{n-1-i},\eeq
where
$$u=\log \frac {\oo_{\varphi}^n}{\oo^n}+\varphi-h_{\oo}.$$By direct calculation, we have
 \beqn \nu_{\oo}(\varphi)-F_{\oo}(\varphi)&=&\frac
1V\int_M\; u\oo^{n}_{\varphi}+\frac
1V\int_M\;h_{\oo}\oo^n+\log\Big(\frac 1V\int_M\;
e^{h_{\oo}-\varphi}\oo^n\Big)\nonumber\\&=&\frac 1V\int_M\;
u\oo^{n}_{\varphi}+\frac 1V\int_M\;h_{\oo}\oo^n+\log\Big(\frac
1V\int_M\; e^{-u}\oo^n_{\varphi}\Big).\label{FK} \eeqn Using
Jensen's inequality, we have
$$\log\Big(\frac
1V\int_M\; e^{-u}\oo^n_{\varphi}\Big)\geq -\frac 1V\int_M\; u\;
\oo_{\varphi}^n.$$ Thus, we have
$$\nu_{\oo}(\varphi)\geq F_{\oo}(\varphi)+\frac
1V\int_M\; h_{\oo}\oo^n.$$
\end{proof}

Now we assume that the $K$-energy is bounded from below. For any
metric $\oo'=\oo+\pbp \varphi_0$, we consider the solution
$\varphi(t)$ of K\"ahler-Ricci flow with the initial metric $\oo'$:
$$\pd {\varphi}{t}=u=\log \frac {\oo_{\varphi}^n}{\oo^n}+\varphi-h_{\oo},\qquad\varphi(0)=\varphi_0.$$
Since $F(t)=F_{\oo}(\varphi(t))$ is decreasing along the
K\"ahler-Ricci flow (cf. \cite{[ChenTian]}), we will prove that
$\nu_{\oo}(t)-F(t)$ has a uniform bound as $t\ri \infty$, and the
bound is independent of the initial metric $\oo'.$ Thus, $F$ is also
bounded from below.\\

Since $F(t)$ is decreasing along the K\"ahler-Ricci
flow, for any $s<t$ by the equality (\ref{FK}) we have \beqn F_{\oo}(\oo')=F(0)&\geq &F(t)-\nu_{\oo}(t)+\nu_{\oo}(t)\nonumber\\
&=&F(t)-\nu_{\oo}(t)+\nu_{\oo}(s)-\frac 1V\int_s^t\int_M\; |\Na u|^2
\oo_{\varphi}^n\nonumber\\
&=&-f(t)+\nu_{\oo}(s)-\frac 1V\int_s^t\int_M\; |\Na u|^2
\oo_{\varphi}^n-\frac 1V\int_M\; h_{\oo}\oo^n.\label{a1}\eeqn where
\beq f(t)=\frac 1V\int_M\; u\;\oo_{\varphi}^n+\log \Big(\frac
1V\int_M\; e^{-u}\oo_{\varphi}^n\Big).\label{eq:f}\eeq If we can
find a sequence of times $t_m\ri \infty$ such that \beq\lim_{m\ri
\infty}f(t_m)=0,\label{a2}\eeq then we can take $t=t_m$ in
(\ref{a1}), and let $m\ri \infty$,
$$F_{\oo}(\oo')\geq \nu_{\oo}(s)-\frac 1V\int_s^{\infty}\int_M\; |\Na u|^2\oo_{\varphi}^n-\frac 1V\int_M\; h_{\oo}\oo^n.$$
Since the $K$-energy is decreasing along K\"ahler-Ricci flow, taking
$s\ri \infty$ in the above inequality we have \beq F_{\oo}(\oo')\geq
\inf\nu_{\oo}-\frac 1V\int_M\; h_{\oo}\oo^n.\label{a3}\eeq Then $F$
is bounded from below. Thus, it suffices to show that
(\ref{a2}) holds.\\

Now we are ready to prove (\ref{a2}).  Since the $K$-energy is
bounded from below, by Lemma \ref{lem2.6} we can normalize the
solution $\varphi(t)$ such that $c(t)>0$ for all $t$, and \beq
\lim_{t\ri\infty}c(t)=\lim_{t\ri\infty}\frac 1V\int_M\;
u\,\oo_{\varphi}^n=0.\label{eq:ct}\eeq By Lemma \ref{Pere}, we prove
\begin{lem}\label{lem3.2}There exists a constant $B$ independent of $t$ such that $|u|\leq B.$
\end{lem}
\begin{proof}We use the notations in Lemma \ref{Pere}. By the
equality (\ref{eq:P1}), we have
$$\int_M\; e^{-u+a_t}\oo_{\varphi}^n=V.$$
It follows that
$$a_t=-\log \Big(\frac 1V\int_M\; e^{-u}\oo_{\varphi}^n\Big).$$
Then Lemma \ref{Pere} implies \beq -A\leq u+\log \Big(\frac
1V\int_M\; e^{-u}\oo_{\varphi}^n\Big)\leq A.\label{eq:3.2}\eeq Since
the $K$-energy is bounded from below, by Lemma \ref{lem2.6} the
integral $\int_M\; u\,\oo_{\varphi}^n$ is uniformly bounded from
above and below. Thus, integrating (\ref{eq:3.2}) we have \beq
\Big|\log \Big(\frac 1V\int_M\; e^{-u}\oo_{\varphi}^n\Big)\Big|\leq
C,\label{eq:3.3}\eeq for some constant $C$. Combining (\ref{eq:3.2})
with (\ref{eq:3.3}), the lemma is proved.
\end{proof}

Next, we prove the following lemma
\begin{lem}\label{lem3.3}For time $t\ri \infty,$
we have $$\un u(t)\ri 0,$$ where $\un u(t)=\frac 1V\int_M\;
u\,e^{h_t}\oo_{\varphi}^n.$ Here we choose $h_t$ as in Lemma
\ref{Pere}.
\end{lem}
\begin{proof}Observe that
\beq \Big(\frac 1V\int_M\; u\;e^{h_t} \oo_{\varphi}^n\Big)^2\leq
\frac 1V\int_M\; u^2e^{h_t} \oo_{\varphi}^n \leq \frac
{e^A}{V}\int_M\; u^2\oo_{\varphi}^n.\label{eq:3.4}\eeq Let
$$b(t)=\int_M\; u^2\oo_{\varphi}^n.$$
Then \beqs \frac{d}{dt}b(t)&=&\int_M\;\Big(2u(\Delta
u+u)+u^2\Delta u\Big)\oo_{\varphi}^n \\
&=&\int_M\;(-2|\Na u|^2+2u^2-2u|\Na u|^2)\oo_{\varphi}^n\\
&\geq &\int_M\;(-2|\Na u|^2+2u^2-u^2-|\Na
u|^4)\oo_{\varphi}^n\\
&\geq &b(t)-(2+A)\int_M\; |\Na u|^2\oo_{\varphi}^n\eeqs where we use
$|\Na u|^2\leq A$ in the last inequality. Thus, integrating the
above inequality from $0$ to $\infty$ we have
$$\int_0^{\infty}\;b(t)dt\leq \limsup_{t\ri \infty}\;
b(t)-b(0)+(2+A)\int_0^{\infty}\;\int_M\; |\Na
u|^2\oo_{\varphi}^n<\infty.
$$ Here the last inequality comes from Lemma \ref{lem3.2} and the fact that the $K$-energy is bounded from below.
By Lemma \ref{Pere}, we have $|\frac d{dt}b(t)|\leq C.$ Hence, we
have $b(t)\ri 0$ as $t\ri \infty$. Therefore, by the inequality
(\ref{eq:3.4}) we have $\un u(t)\ri 0.$
\end{proof}

Now we can prove
\begin{lem}\label{lem3.4}There is a
sequence of times $t_m\ri \infty$ such that $$\lim_{m\ri
\infty}f(t_m)\ri 0,$$ where $f$ is defined by (\ref{eq:f}).
\end{lem}
\begin{proof}By the equalities (\ref{eq:f}) and (\ref{eq:ct}), it
suffices to find a sequence of times $t_m$ such that \beq\lim_{m\ri
\infty}\log\Big(\frac 1V \int_M\;
e^{-u}\oo_{\varphi}^n\Big)\Big|_{t=t_m}= 0\label{eq:3.5}\eeq

Since $u$ and $\un u$ are bounded by Lemma \ref{lem3.2} and Lemma
\ref{lem3.3}, we have the Taylor expansion \beq e^{-(u-\un
u)}=1+\sum_{k=1}^{\infty}\frac {(-1)^k}{k!}(u-\un u)^k.\eeq
Therefore, \beq \log\Big(\frac 1V \int_M\;
e^{-u}\oo_{\varphi}^n\Big)=-\un u+\log
\Big(1+\frac{1}V\int_M\;\sum_{k=1}^{\infty}\frac {(-1)^k}{k!}(u-\un
u)^k\oo_{\varphi}^n \Big).\label{eq:3.6}\eeq
 Now by Lemma \ref{Pere}, we have
\beqn \Big|\int_M\; \sum_{k=1}^{\infty}\frac {(-1)^k}{k!}(u-\un
u)^k\oo_{\varphi}^n\Big| &\leq &\int_M\; \sum_{k=1}^{\infty}\frac
{1}{k!}|u-\un
u|^k\oo_{\varphi}^n\nonumber\\
&\leq &\sum_{k=1}^{\infty}\frac {e^A}{k!}\int_M\; |u-\un
u|^ke^{h_t}\oo_{\varphi}^n. \label{eq:3.7}\eeqn Then by the
Poincar\'e inequality in Lemma \ref{TZ}, we know \beqn
&&\sum_{k=1}^{\infty}\frac {e^A}{k!}\int_M\; |u-\un
u|^ke^{h_t}\oo_{\varphi}^n\nonumber\\ &\leq &e^A\int_M\; |u-\un
u|e^{h_t}\oo_{\varphi}^n+\sum_{k=2}^{\infty}\frac
{e^A(2B)^{k-2}}{k!}\int_M\; |u-\un
u|^2e^{h_t}\oo_{\varphi}^n\nonumber\\
&\leq &e^A\sqrt{V}\sqrt{\int_M\; |u-\un
u|^2e^{h_t}\oo_{\varphi}^n}+\frac {e^{A+2B}}{(2B)^2}\int_M\; | u-\un
u|^2e^{h_t}\oo_{\varphi}^n\nonumber\\&\leq&
e^A\sqrt{V}\sqrt{\int_M\; |\Na u|^2e^{h_t}\oo_{\varphi}^n}+\frac
{e^{A+2B}}{(2B)^2}\int_M\; |\Na u|^2e^{h_t}\oo_{\varphi}^n\nonumber\\
&\leq&e^{\frac 32A}\sqrt{V}\sqrt{\int_M\; |\Na
u|^2\oo_{\varphi}^n}+\frac {e^{2A+2B}}{(2B)^2}\int_M\; |\Na
u|^2\oo_{\varphi}^n.\label{eq:3.8} \eeqn Since the $K$-energy is
bounded from below, by (\ref{eq:dK}) we can find a sequence of times
$t_m\ri \infty$ such that
$$\int_M\;|\Na u|^2\oo_{\varphi}^n\Big|_{t=t_m}\ri 0.$$
Combining this with (\ref{eq:3.6})-(\ref{eq:3.8}), we know
(\ref{eq:3.5}) holds. The lemma is proved.
\end{proof}

By Lemma \ref{lem3.4}, the equality (\ref{a2}) holds. This implies
$F_{\oo}$ is bounded from below and the inequality (\ref{a3}) holds.
Combining this with Lemma \ref{lem3.1}, the main theorem is proved.

School of Mathematical Sciences, Peking University, Beijing, 100871,
P. R. China\\lihaozhao@gmail.com


\begin{thebibliography}{2}
\bibitem{[BaMa]}S. Bando, T. Mabuchi: {\it Uniqueness of Einstein K\"ahler metrics modulo connected group actions},  Algebraic geometry, Sendai, 1985,  11--40, Adv. Stud. Pure Math., 10, North-Holland, Amsterdam, 1987.
\bibitem{[Chen0]}X. X. Chen: {\it On the lower bound of the Mabuchi
energy and its applications},  Int. Math. Res. Lett. 12 (2000)
607-623.

\bibitem{[Chen1]}X. X. Chen: {\it On the lower bound of energy functional $E_1 (I)$--
 a stability theorem on the K\"ahler-Ricci flow}, J. Geometric
 Analysis. 16 (2006) 23-38.
\bibitem{[CLW]}X.X. Chen, H. Li, B. Wang: {\it On the K\"ahler-Ricci flow with small initial $E_1$ energy
(I)}, math.DG/0609694.

\bibitem{[ChenTian]}X. X. Chen, G. Tian: {\it Ricci flow on K\"ahler-Einstein surfaces},
  Invent. Math.  147  (2002),  no. 3, 487--544.
\bibitem{[chen-tian1]}X. X. Chen, G. Tian: {\it Ricci flow on K\"ahler-Einstein
manifolds}, Duke. Math. J.  131, (2006),  no. 1, 17-73.


\bibitem{[DingTian]}W. Y. Ding, G. Tian: {\it The generalized Moser-Trudinger inequality},
 Proceedings of Nankai International Conference of Nonlinear Analysis, 1993.
\bibitem{[Futaki]}A. Futaki: {\it K\"ahler-Einstein metrics and integral invariants.} Lecture Notes in Math. 1314, Springer 1988.

\bibitem{[Li]}H. Li: {\it A new formula for the energy functionals $E_k$ and its
applications}, math.DG/0609724.
\bibitem{[Ma]}T. Mabuchi: {\it $K$-energy maps integrating Futaki invariants},
 Tohoku Math. J. (2) 38(1986), no. 4, 575-593.
\bibitem{[NaTi]}N. Sesum, G. Tian: {\it Bounding scalar curvature and diameter along the
K\"ahler-Ricci flow (after Perelman) and some applications},
preprint.

\bibitem{[Pali]}N. Pali. {\it A consequence of a lower bound of the $K$-energy},  Int. Math. Res. Not.  2005,  no. 50, 3081--3090.

\bibitem{[Natasa]}G.  Perelman: {\it Unpublished work on K\"ahler-Ricci flow}.

\bibitem{[R]}Y. Rubinstein: {\it On energy functionals and the existence of K\"ahler-Einstein
metrics}. math.DG/0612440.

\bibitem{[SoWe]}J. Song, B. Weinkove: {\it Energy functionals and canonical K\"ahler metrics}, math.DG/0505476.
\bibitem{[Tian97]}G. Tian: {\it K\"ahler-Einstein metrics with positive scalar curvature},  Invent. Math.  130  (1997),  no. 1, 1--37.

\bibitem{[Tianbook]}G. Tian:  Canonical metrics in K\"ahler geometry. Notes taken by Meike Akveld. Lectures in Mathematics ETH Z\"urich. Birkh\"auser Verlag, Basel, 2000.
\bibitem{[TianZhu]}G. Tian, X.H. Zhu: {\it Convergence of K\"ahler-Ricci
flow}, J. Amer. Math. Soc. 20 (2007), no. 3, 675-699.
\bibitem{[Tosatti]}V. Tosatti: {\it On the Critical Points of the $E_k$ Functionals in K\"ahler
Geometry}, math.DG/0506021.




\end{thebibliography}
\end{document}